\documentclass[12pt]{article}

\usepackage[dvipdfm, bookmarksopen=true]{hyperref}
\usepackage[centertags]{amsmath}
\usepackage{amsfonts}
\usepackage{amssymb}
\usepackage{amsthm}

\theoremstyle{plain}
\newtheorem{thm}{Theorem}[section]
\newtheorem{cor}[thm]{Corollary}
\newtheorem{lem}[thm]{Lemma}

\newtheorem{prop}[thm]{Proposition}

\theoremstyle{definition}

\theoremstyle{remark}
\newcommand{\beast}{\begin{eqnarray*}}
\newcommand{\eeast}{\end{eqnarray*}}

\title{{\bf Two-lit trees for lit-only sigma-game}
}

\author{Hau-wen Huang}

\date{\today}

\begin{document}
\maketitle

\begin{abstract}
A configuration of the lit-only $\sigma$-game on a finite graph $\Gamma$ is an assignment of one of two states, {\it on} or {\it off}, to all vertices of $\Gamma.$ Given a configuration, a move of the lit-only $\sigma$-game on $\Gamma$ allows the player to choose an {\it on} vertex $s$ of $\Gamma$ and change the states of all neighbors of $s.$
Given any integer $k$, we say that $\Gamma$ is $k$-lit if, for any configuration, the number of {\it on} vertices can be reduced to at most $k$ by a finite sequence of moves. Assume that $\Gamma$ is a tree with a perfect matching. We show that $\Gamma$ is $1$-lit and any tree obtained from $\Gamma$ by adding a new vertex on an edge of $\Gamma$ is $2$-lit.
\end{abstract}

{\footnotesize{\bf Keywords:} group action, lit-only sigma-game, symplectic forms}

{\footnotesize {\bf 2010  MSC Primary:} 05C57; {\bf Secondary:} 15A63, 20F55}\\

\section{Introduction}\label{introduction}

In 1989, Sutner \cite{ks:89} introduced a one-player game called the $\sigma$-game. The $\sigma$-game is played on a finite directed graph $\Gamma$ without multiple edges. A {\it configuration} of the $\sigma$-game on $\Gamma$ is an assignment of one of two states, {\it on} or {\it off}, to all vertices of $\Gamma.$  Given a configuration, a {\it move} of the $\sigma$-game on $\Gamma$ allows the player to pick any vertex $s$ of $\Gamma$ and change the states of all neighbors of $s.$ Given an initial configuration, the goal is to minimize the number of {\it on} vertices of $\Gamma$ or to reach an assigned configuration by a finite sequence of moves. If only {\it on} vertex can be chosen in each move, we come to the variation: {\it lit-only $\sigma$-game}. The goal of the lit-only $\sigma$-game is the same as that of the $\sigma$-game. Given an integer $k$, we say that $\Gamma$ is {\it $k$-lit} if, for any configuration, the number of {\it on} vertices can be reduced to at most $k$ by a finite sequence of moves of the lit-only $\sigma$-game on $\Gamma$. Motivated by the goal of the lit-only $\sigma$-game, we are interested in the smallest integer $k$, the {\it minimum light number} of $\Gamma$ \cite{xwykw:07}, for which $\Gamma$ is $k$-lit.

As far as we know, the notion of the lit-only $\sigma$-game first implicitly occurred in the classification of the equivalence classes of Vogan diagrams, which implies that all simply-laced Dynkin diagrams, the trees shown as below, are $1$-lit (cf. \cite{abjs:49,mkc:04}). Extending this result, Wang and Wu proved that any tree with $k$ leaves is $\lceil k/2\rceil$-lit (cf. \cite[Theorem~3]{xwykw:07}). Their recent result gave an insight into the  difference between the $\sigma$-game and the lit-only $\sigma$-game on trees with zero or more loops (cf. \cite[Theorem~14]{xwykw:09}). As a consequence, the trees with perfect matchings are $2$-lit. The first main result of this paper improves this consequence.

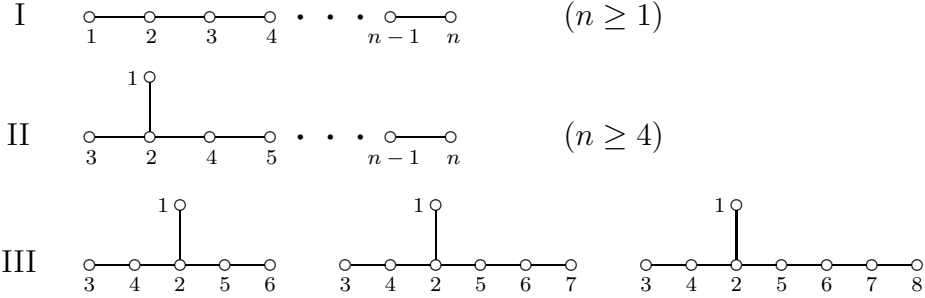
\setlength{\unitlength}{1mm}
\begin{picture}(150,50)

\put(17,38){I}
\multiput(27,39)(8,0){4}{\circle{1.5}}
\put(26.5,35.5){{\scriptsize $1$}}
\put(34.5,35.5){{\scriptsize $2$}}
\put(42.5,35.5){{\scriptsize $3$}}
\put(50.5,35.5){{\scriptsize $4$}}
\multiput(55,39)(4,0){3}{\circle*{0.8}}
\multiput(67,39)(8,0){2}{\circle{1.5}}
\put(64,35.5){{\scriptsize $n-1$}}
\put(74.5,35.5){{\scriptsize $n$}}
\multiput(27.75,39)(8,0){3}{\line( 1, 0){6.5}}
\put(67.75,39){\line( 1, 0){6.5}}
\put(90,38){$(n\geq 1)$}

\put(16,22){II}
\multiput(27,23)(8,0){4}{\circle{1.5}}
\put(26.5,19.5){{\scriptsize $3$}}
\put(34.5,19.5){{\scriptsize $2$}}
\put(42.5,19.5){{\scriptsize $4$}}
\put(50.5,19.5){{\scriptsize $5$}}
\multiput(55,23)(4,0){3}{\circle*{0.8}}
\put(35,31){\circle{1.5}}
\put(32,30){{\scriptsize $1$}}
\multiput(67,23)(8,0){2}{\circle{1.5}}
\put(64,19.5){{\scriptsize $n-1$}}
\put(74.5,19.5){{\scriptsize $n$}}
\multiput(27.75,23)(8,0){3}{\line( 1, 0){6.5}}
\put(35,23.75){\line( 0, 1){6.5}}
\put(67.75,23){\line( 1, 0){6.5}}
\put(90,22){$(n\geq 4)$}

\put(15.25,5){III}
\multiput(27,6)(6,0){5}{\circle{1.5}}
\put(26.25,2.5){{\scriptsize $3$}}
\put(39,14){\circle{1.5}}
\put(36,13){{\scriptsize $1$}}
\multiput(27.75,6)(6,0){4}{\line( 1, 0){4.5}}
\put(32.25,2.5){{\scriptsize $4$}}
\put(38.25,2.5){{\scriptsize $2$}}
\put(44.25,2.5){{\scriptsize $5$}}
\put(50.25,2.5){{\scriptsize $6$}}
\put(39,6.75){\line( 0, 1){6.5}}

\multiput(61,6)(6,0){6}{\circle{1.5}}
\put(60.25,2.5){{\scriptsize $3$}}
\put(73,14){\circle{1.5}}
\put(70,13){{\scriptsize $1$}}
\multiput(61.75,6)(6,0){5}{\line( 1, 0){4.5}}
\put(66.25,2.5){{\scriptsize $4$}}
\put(72.25,2.5){{\scriptsize $2$}}
\put(78.25,2.5){{\scriptsize $5$}}
\put(84.25,2.5){{\scriptsize $6$}}
\put(90.25,2.5){{\scriptsize $7$}}
\put(73,6.75){\line( 0, 1){6.5}}

\multiput(101,6)(6,0){7}{\circle{1.5}}
\put(100.25,2.5){{\scriptsize $3$}}
\put(113,14){\circle{1.5}}
\put(110,13){{\scriptsize $1$}}
\multiput(101.75,6)(6,0){6}{\line( 1, 0){4.5}}
\put(106.25,2.5){{\scriptsize $4$}}
\put(112.25,2.5){{\scriptsize $2$}}
\put(118.25,2.5){{\scriptsize $5$}}
\put(124.25,2.5){{\scriptsize $6$}}
\put(130.25,2.5){{\scriptsize $7$}}
\put(136.25,2.5){{\scriptsize $8$}}
\put(113,6.75){\line( 0, 1){6.5}}

\end{picture}

The lit-only $\sigma$-game on a finite simple graph $\Gamma$ can be regarded as a representation of the simply-laced Coxeter group associated with $\Gamma$ (cf. \cite{hw:08-1}). From this viewpoint, we apply some results from \cite{rd:05} to show that the trees with perfect matchings, except the paths of even order, are $1$-lit. Combining this with the result that all paths, namely the trees in class I, are $1$-lit, our first result can be simply stated as follows.

\begin{thm}\label{thm1.2}
Any tree with a perfect matching is $1$-lit.
\end{thm}

Theorem~\ref{thm1.2} gives a large family of $1$-lit trees containing the first and third trees in class III. It is natural to ask if there is also a large family of $1$-lit trees containing the second tree in class III. This question motivates the discovery of our second result. Assume that $\Gamma$ is a tree with a perfect matching $\mathcal{P}.$ An {\it alternating path} in $\Gamma$ {\rm (}with respect to $\mathcal{P}${\rm )} is a path in which the edges belong alternatively to $\mathcal{P}$ and
not to $\mathcal{P}.$ For each vertex $s$ of $\Gamma$, we define $a_s$ to be the number of the alternating paths starting from the edge in $\mathcal P$ incident to $s$ and ending on some edge in $\mathcal{P}.$ An edge of $\Gamma$ is said to be of {\it odd  {\rm(resp.} {\it even}{\rm)}  type} if its two endpoints $s,t$ satisfy that $a_s+a_t$ is odd (resp. even). We make use of algebraic and linear algebraic techniques to show that

\begin{thm}\label{thm1.4}
Assume that $\Gamma$ is a tree with a perfect matching. Then  the tree obtained from $\Gamma$ by inserting a vertex on an edge of odd {\rm(}resp.  even{\rm)} type is $1$-lit  {\rm(}resp. $2$-lit{\rm)}.
\end{thm}

Let $\Gamma$ denote the first tree in class III. The edge of $\Gamma$ joining $5$ and $6$ is of odd type because of $a_5=1$ and $a_6=2$. Therefore the second tree in class III does be a special case of Theorem~\ref{thm1.4}. On the other hand, if we add a vertex on the edge of $\Gamma$ between $1$ and $2$, the resulting tree is not $1$-lit by \cite[Proposition 3.2]{mkc:06}. Therefore, in general, for any tree $\Gamma$ with a perfect matching, the tree obtained from $\Gamma$ by adding a vertex on an edge is not $1$-lit.

The statements of Theorem~\ref{thm1.2} and Theorem~\ref{thm1.4} are combinatorial. It is reasonable to believe that these results can be proved by combinatorial arguments. In addition, motivated by Theorem~\ref{thm1.4}, we would like to ask if given a tree $\Gamma$ with a perfect matching, any subdivision of $\Gamma$ is $2$-lit. We leave these as open problems.

\section{Preliminaries}
For the rest of this paper, let $\Gamma=(S,R)$ denote a finite simple graph with vertex set $S$ and edge set $R$. The edge set $R$ is a set of some $2$-element subsets of $S$. For any distinct $s,t\in S$, we write $st$ or $ts$ to denote the $2$-element subset $\{s,t\}$ of $S$. Let $\mathbb{F}_2$ denote the two-element field $\{0,1\}$. Let $V$ denote a $\mathbb{F}_2$-vector space that has a basis $\{\alpha_s~|~s\in S\}$ in one-to-one correspondence with $S.$ Let $V^*$ denote the dual space of $V$. For each $s\in S$, define $f_s\in V^*$ by
\begin{align}\label{e1.1}
f_s(\alpha_t)=\left\{
\begin{array}{ll}
1 \qquad &\hbox{if $s=t,$}\\
0 \qquad &\hbox{if $s\not=t$}
\end{array}
\right.
\end{align}
 for all $t\in S$. The set $\{f_s~|~s\in S\}$ is a basis of $V^*$ and called the basis of $V^*$ dual to $\{\alpha_s~|~s\in S\}$. Each configuration $f$ of the lit-only $\sigma$-game on $\Gamma$ is interpreted as the vector
\begin{align}\label{e1.4}
\sum_{\scriptsize \hbox{ {\it on} vertices $s$}} f_s\in V^*,
\end{align}
if all vertices of $\Gamma$ are assigned the {\it off} state by $f,$ we interpret (\ref{e1.4}) as the zero vector of $V^*.$ For any $s\in S$ and $f\in V^*,$ $f(\alpha_s)=1$ (resp. 0) means that the vertex $s$ is assigned the {\it on} (resp. {\it off} ) state by $f.$ For each $s\in S$ define a linear transformation $\kappa_s:V^*\rightarrow V^*$ by
\begin{align}\label{e1.3}
\kappa_sf= f+f(\alpha_s)\sum_{st\in R}f_{t} \qquad \quad &\hbox{for all $f\in V^*.$}
\end{align}
Fix a vertex $s$ of $\Gamma.$ Given any $f\in V^*,$ if the state of $s$ is {\it on} then $\kappa_sf$ is obtained from $f$ by changing the states of all neighbors of $s;$ if the state of $s$ is {\it off} then $\kappa_sf=f.$ Therefore we may view $\kappa_s$ as the move of the lit-only $\sigma$-game on $\Gamma$ for which we choose the vertex $s$ and change the states of all neighbors of $s$ if the state of $s$ is {\it on}. In particular $\kappa_s^2=1$, the identity map on $V^*$, and so $\kappa_s\in {\rm GL}(V^*)$, the general linear group of $V^*$.

The simply-laced Coxeter group $W$ associated with $\Gamma=(S,R)$ is a group generated by the set $S$ subject to the following relations:
\begin{align}
s^2&=1,\label{defnW1}\\
(st)^2&=1 \qquad \hbox{if $st\not\in R,$} \label{defnW2}\\
(st)^3&=1 \qquad \hbox{if $st\in R$} \label{defnW3}
\end{align}
for all $s,t\in S$. By \cite[Theorem~3.2]{hw:08-1}, there is a unique representation $\kappa:W\to {\rm GL}(V^*)$ such that $\kappa(s)=\kappa_s$ for all $s\in S$. For any $f,g\in V^*,$ observe that $g$ can be obtained from $f$ by a finite sequence of moves of the lit-only $\sigma$-game on $\Gamma$ if and only if there exists $w\in W$ such that $g=\kappa(w)f.$ In view of this we define an action of $W$ on $V^*$ by
$$
wf=\kappa(w)f \qquad \quad \hbox{ for all $w\in W$ and $f\in V^*.$}
$$
In terms of our terminology, given an integer $k$, the simple graph $\Gamma$ is $k$-lit if and only if for any $W$-orbit $O$ of $V^*$, there exists a subset $K$ of $S$ with cardinality at most $k$ such that $\sum_{s\in K} f_s\in O.$

Let $B:V\times V\to \mathbb{F}_2$ denote the symplectic form defined by
\begin{align}\label{e1.2}
\begin{array}{ll}
B(\alpha_s,\alpha_t)=\left\{
\begin{array}{ll}
1 \qquad &\hbox{if $st\in R,$}\\
0 \qquad &\hbox{else}
\end{array}
\right.
\end{array}
\end{align}
for all $s,t\in S$. By (\ref{e1.1}) and (\ref{e1.2}), for all $s\in S$ and $\alpha \in V$ we have
\begin{align}\label{e2.6}
B(\alpha_s,\alpha)=\sum_{st\in R}f_t(\alpha).
\end{align}
The {\it radical} of $V$ (relative to $B$), denoted by ${\rm rad}\hspace{0.05cm}V$, is the subspace of $V$ consisting of the $\alpha\in V$ that satisfy $B(\alpha,\beta)=0$ for all $\beta\in V$. The form $B$ is said to be {\it degenerate} if ${\rm rad}\hspace{0.05cm}V\not=\{0\}$ and {\it nondegenerate} otherwise. The graph $\Gamma$ is said to be {\it degenerate} (resp. {\it nondegenerate}) if $B$ is degenerate (resp. nondegenerate). The form $B$ induces a linear map $\theta:V\rightarrow V^*$ given by
\begin{align}\label{e4.4}
\theta(\alpha)\beta=B(\alpha,\beta)  \qquad \quad &\hbox{for all $\alpha,\beta\in V.$}
\end{align}
By (\ref{e2.6}) and (\ref{e4.4}), for each $s\in S$ we have
\begin{align}\label{e4.6}
\theta(\alpha_s)=\sum_{st\in R}f_t.
\end{align}
Let $A$ denote the adjacency matrix of $\Gamma$ over $\mathbb{F}_2$. Observe that the kernel of $\theta$ is ${\rm rad}\hspace{0.05cm}V$ and the matrix representing $B$ with respect to $\{\alpha_s~|~s\in S\}$ is exactly $A$. Therefore we have

\begin{lem}\label{lem4.2}
The following statements are equivalent:
\begin{enumerate}
\item $\Gamma$ is a nondegenerate graph.

\item $\theta$ is an isomorphism of vector spaces.

\item $A$ is invertible.
\end{enumerate}
\end{lem}

The determinant of $A$ is $0$ (resp. $1$) if and only if the number of perfect matchings in $\Gamma$ is even (resp. odd) (see \cite[Section 2.1]{cdg:93} for example). Combining this with Lemma~\ref{lem4.2} we have

\begin{prop}\label{lem3.1}
The following statements are equivalent:
\begin{enumerate}
\item $\Gamma$ is a nondegenerate graph.

\item The number of perfect matchings in $\Gamma$ is odd.
\end{enumerate}
\end{prop}

Since a tree contains at most one perfect matching and by Proposition~\ref{lem3.1}, we have

\begin{cor}\label{prop3.1}
The following statements are equivalent:
\begin{enumerate}
\item $\Gamma$ is a nondegenerate tree.

\item $\Gamma$ is a tree with a perfect matching.
\end{enumerate}
\end{cor}

\begin{prop}\label{prop3.2}
{\rm (\cite[Lemma~2.4]{eigen_nondegtree:02}).} Assume that $\Gamma$ is a tree of order at least four and with a perfect matching. Then there exist two vertices of $\Gamma$ with degree two.
\end{prop}

\section{Proof of Theorem~\ref{thm1.2}}
The lit-only $\sigma$-game is closely related to another combinatorial game. We call this game the {\it Reeder's game} because as far as we know, this game first appeared in one of Reeder's papers \cite{rd:05}. The Reeder's game is a one-player game played on a finite simple graph $\Gamma$. A configuration of the Reeder's game on $\Gamma$ is an assignment of one of two states, {\it on} or {\it off}, to each vertex of $\Gamma$. Given a configuration, a move of the Reeder's game on $\Gamma$ consists of choosing a vertex $s$ and changing the state of $s$ if the number of {\it on} neighbors of $s$ is odd. Given an initial configuration, the goal is to minimize the number of {\it on} vertices of $\Gamma$ by a finite sequence of moves of the Reeder's game on $\Gamma.$

We interpret each configuration $\alpha$ of the Reeder's game on $\Gamma$ as the vector
\begin{align}\label{e2.5}
\sum_{\scriptsize \hbox{{\it on} vertices $s$}}\alpha_s\in V,
\end{align}
if all vertices of $\Gamma$ are assigned the {\it off} state by $\alpha,$ we interpret (\ref{e2.5}) as the zero vector of $V.$ For any $\alpha\in V,$ observe that $f_s(\alpha)=1$ (resp. $0$) means that the vertex $s$ is assigned the {\it on} (resp. {\it off}) state by $\alpha.$ For each $s\in S$ define a linear transformation $\tau_s:V\rightarrow V$ by
\begin{align}\label{e2.3}
\tau_s\alpha=\alpha+B(\alpha_s,\alpha)\alpha_s \qquad \quad &\hbox{for all $\alpha\in V.$}
\end{align}
Fix a vertex $s$ of $\Gamma.$ By (\ref{e2.6}), for any $\alpha\in V,$ if the number of {\it on} neighbors of $s$ is odd then $\tau_s\alpha$ is obtained from $\alpha$ by changing the state of $s;$ if the number of {\it on} neighbors of $s$ is even then $\tau_s\alpha=\alpha.$ Therefore we may view $\tau_s$ as the move of the Reeder's game on $\Gamma$ for which we choose the vertex $s$ and change the state of $s$ if the number of {\it on} neighbors of $s$ is odd. In particular $\tau_s^2=1$, the identity map on $V$, and so $\tau_s\in {\rm GL}(V)$, the general linear group of $V$.

By \cite[Section 5]{rd:05}, there exists a unique representation $\tau:W\rightarrow {\rm GL}(V)$ such that $\tau(s)=\tau_s$ for all $s\in S.$ For any $\alpha,\beta\in V,$ observe that $\beta$ can be obtained from $\alpha$ by a finite sequence of moves of the Reeder's game on $\Gamma$ if and only if there exists $w\in W$ such that $\beta=\tau(w)\alpha.$ In view of this we define an action $W$ on $V$ by
$$
w\alpha=\tau(w)\alpha \qquad \quad \hbox{for all $w\in W$ and $\alpha\in V.$}
$$
A quadratic form $Q:V\rightarrow \mathbb{F}_2$, given in \cite[Section~1]{rd:05},  is defined by
\begin{eqnarray}
Q(\alpha_s)&=&1  \hspace{5.9cm}\hbox{for all $s\in S,$} \label{e2.1}\\
Q(\alpha+\beta)&=&Q(\alpha)+Q(\beta)+B(\alpha,\beta)  \hspace{1.9cm} \hbox{for all $\alpha,\beta\in V.$} \label{e2.2}
\end{eqnarray}
Observe that $\tau$ preserves $Q$, namely
\begin{align}\label{e2.7}
Q(\tau(w)\alpha)=Q(\alpha) \qquad \quad \hbox{for all $w\in W$ and $\alpha\in V$}.
\end{align}
The {\it kernel} of $Q$, denoted by $ {\rm Ker}\hspace{0.05cm}Q$, is the subspace of ${\rm rad}\hspace{0.05cm}V$ consisting of all $\alpha\in {\rm rad}\hspace{0.05cm}V$ that satisfy $Q(\alpha)=0$. The {\it orthogonal group} $O(V)$ (relative to $Q$) is the subgroup of ${\rm GL}(V)$ consisting of the $\sigma\in {\rm GL}(V)$ such that $Q(\sigma\alpha)=Q(\alpha)$ for all $\alpha\in V.$

\begin{lem}\label{thm2.2} {\rm (\cite[Section~2; Theorem 7.3]{rd:05}).}
Let $\Gamma$ denote a tree which is not a path. Assume that $ {\rm Ker}\hspace{0.05cm}Q$ is equal to $\{0\}$. Then
$
\tau(W)=O(V).
$
Moreover the $W$-orbits on $V$ are
\begin{gather*}
Q^{-1} (1)\setminus {\rm rad}\hspace{0.5mm}V, \qquad \quad Q^{-1}(0)\setminus \{0\}, \qquad \quad  \{\alpha\} \quad \hbox{for all $\alpha\in {\rm rad}\hspace{0.5mm}V$}.
\end{gather*}
\end{lem}

As a corollary of Lemma~\ref{thm2.2} we have

\begin{cor}\label{cor2.1}
Assume that $\Gamma$ is a nondegenerate tree which is not a path. Then the $W$-orbits on $V$ are
\begin{gather*}
Q^{-1} (1), \qquad \quad Q^{-1}(0)\setminus \{0\}, \qquad \quad  \{0\}.
\end{gather*}
\end{cor}

Recall that the transpose of a linear transformation $\sigma:V\to V$ is the linear transformation ${}^t\sigma:V^*\to V^*$ defined by $({}^t\sigma f)(\alpha)=f(\sigma \alpha)$ for all $f\in V^*$ and $\alpha\in V$.

\begin{lem}\label{prop2.1}
The representation $\kappa$ is the dual representation of $\tau.$
\end{lem}
\begin{proof}
Let $s\in S$ be given. Using (\ref{e1.3}), (\ref{e2.6}) and (\ref{e2.3}),  we find that $(\kappa_s f)(\alpha)=({}^t\tau_sf)(\alpha)$ for all $f\in V^*$ and $\alpha\in V$. Therefore $\kappa_s={}^t\tau_s$. Since the elements $s\in S$ generate $W$ and $s^{-1}=s$ in $W$, we have $\kappa(w)={}^t\tau(w^{-1})$ for all $w\in W$. The result follows.
\end{proof}

\begin{lem}\label{lem2.3}
For all $w\in W$ and $\alpha,\beta\in V$ we have
\begin{gather*}
B(\tau(w)\alpha,\tau(w)\beta)=B(\alpha,\beta).
\end{gather*}
\end{lem}
\begin{proof}
Fix $s\in S$. Pick any $\alpha,\beta\in V$. Using (\ref{e1.2}), (\ref{e2.3}) to simplify $B(\tau_s\alpha,\tau_s\beta)$ we obtain that $B(\tau_s\alpha,\tau_s\beta)=B(\alpha,\beta)$. The result follows since the elements $s\in S$ generate $W$.
\end{proof}

We have seen that $\kappa$ is the dual representation of $\tau$ and that $\tau$ preserves the form $B$. By the principles of representation theory, the following lemma is straightforward. For the convenience of the reader we include the proof.

\begin{lem}\label{lem4.3}
$\kappa(w)\circ\theta=\theta\circ\tau(w)$ for all $w\in W.$
\end{lem}
\begin{proof}
Let $w\in W$ be given. Replacing $\beta$ by $\tau(w^{-1})\beta$ in Lemma~\ref{lem2.3}, we obtain
\begin{align}\label{e4.5}
B(\tau(w)\alpha,\beta)=B(\alpha,\tau(w^{-1})\beta) \qquad \quad \hbox{for all $\alpha,\beta\in V.$}
\end{align}
Using (\ref{e4.4}) we can rewrite (\ref{e4.5}) as
\begin{align}\label{e4.1}
(\theta\circ \tau(w))(\alpha)=({}^t\tau(w^{-1})\circ\theta)(\alpha) \qquad \quad \hbox{for all $\alpha\in V.$}
\end{align}
By Lemma~\ref{prop2.1} the right-hand side of (\ref{e4.1}) is equal to $(\kappa(w)\circ \theta)(\alpha).$ The result follows.
\end{proof}

As a consequence of Lemma~\ref{lem4.3} we have

\begin{cor}\label{lem4.4}
Assume that $\theta$ is an isomorphism of vector spaces. Then the representation $\tau$ is equivalent to the representation $\kappa$ via $\theta.$
Moreover the map from the $W$-orbits of $V$ to the $W$-orbits of $V^*$ defined by
\begin{align*}
O\mapsto \theta(O) \qquad \quad &\hbox{for all $W$-orbits $O$ of $V$}
\end{align*}
is a bijection.
\end{cor}

Combining Lemma~\ref{lem4.2}, Corollary~\ref{cor2.1} and Corollary~\ref{lem4.4}, we have

\begin{cor}\label{cor4.1}
Assume that $\Gamma$ is a nondegenerate tree which is not a path. Then the $W$-orbits of $V^*$ are
\begin{align*}
\theta(Q^{-1}(1)), \qquad \quad  \theta(Q^{-1}(0))\setminus\{0\}, \qquad \quad \{0\}.
\end{align*}
\end{cor}

Our last tool for proving Theorem~\ref{thm1.2} is \cite[Theorem 6]{gxw:08}. Here we offer a short proof of this result.

\begin{lem}\label{thm4.1}{\rm (\cite[Theorem 6]{gxw:08}).}
Assume that $\Gamma=(S,R)$ is a nondegenerate graph. Let $s\in S$ and let $f\in V^*$ with $f(\alpha_s)=0.$ Then $f$ and $f+\sum_{st\in R}f_t$ are in distinct $W$-orbits of $V^*.$
\end{lem}
\begin{proof}
Suppose on the contrary that there exists $w\in W$ such that
\begin{align}\label{e4.2}
\kappa(w)f=f+\sum_{st\in R}f_t.
\end{align}
Since $\theta$ is a bijection by Lemma~\ref{lem4.2}, there exists a unique $\alpha\in V$ such that $\theta(\alpha)=f.$ By (\ref{e4.6}), we can rewrite  (\ref{e4.2}) as $\kappa(w)(\theta(\alpha))=\theta(\alpha+\alpha_s)$. By Lemma~\ref{lem4.3} and since $\theta$ is a bijection, we obtain
\begin{align}\label{e4.3}
\tau(w)\alpha=\alpha+\alpha_s.
\end{align}
We now apply $Q$ to either side of (\ref{e4.3}).
By (\ref{e2.7}), the left-hand side is equal to $Q(\alpha).$ By (\ref{e4.4}) and the assumption on $f$, we have $B(\alpha,\alpha_s)=0.$ By this and using (\ref{e2.1}) and (\ref{e2.2}), we find that the right-hand side is equal to $Q(\alpha)+1,$ a contradiction.
\end{proof}

It is now a simple matter to prove Theorem~\ref{thm1.2}.

\medskip

\noindent {\bf Proof of Theorem~\ref{thm1.2}:} Let $\Gamma$ be a tree with a perfect matching. Recall from Section~\ref{introduction} that all paths are $1$-lit. Thus it is enough to treat the case that $\Gamma$ is not a path. Such a $\Gamma$ has order at least four. By Proposition~\ref{prop3.2} there exists a vertex $s$ of $\Gamma$ with degree two. Let $u,$ $v$ denote the neighbors of $s.$ By Corollary~\ref{prop3.1} the tree $\Gamma$ is nondegenerate. Applying Lemma~\ref{thm4.1} to $f=f_u,$ we obtain $f_u$ and $f_v$ in distinct $W$-orbits of $V^*.$ Since there are exactly two nonzero $W$-orbits of $V^*$ by Corollary~\ref{cor4.1}, this implies that $\Gamma$ is $1$-lit. \hfill $\Box$

\medskip

The following example gives a nondegenerate graph which is not $1$-lit. Let $\Gamma=(S,R)$ be the graph shown as follows.

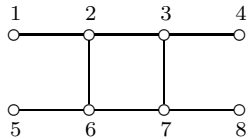
\setlength{\unitlength}{1mm}
\begin{picture}(150,25)

\put(59.5,19){{\scriptsize $1$}}
\put(69.5,19){{\scriptsize $2$}}
\put(79.5,19){{\scriptsize $3$}}
\put(89.5,19){{\scriptsize $4$}}
\multiput(60,17)(10,0){4}{\circle{1.5}}
\put(59.5,3.5){{\scriptsize $5$}}
\put(69.5,3.5){{\scriptsize $6$}}
\put(79.5,3.5){{\scriptsize $7$}}
\put(89.5,3.5){{\scriptsize $8$}}
\multiput(60,7)(10,0){4}{\circle{1.5}}
\multiput(60.75,17)(10,0){3}{\line( 1, 0){8.5}}
\multiput(60.75,7)(10,0){3}{\line( 1, 0){8.5}}
\multiput(70,7.75)(10,0){2}{\line( 0, 1){8.5}}
\end{picture}

\noindent  The graph $\Gamma$ contains the only perfect matching $\{\{1,2\},\{3,4\},\{5,6\},\{7,8\}\}$. By Proposition~\ref{lem3.1} the graph $\Gamma$ is nondegenerate. Let $f=f_2+f_3+f_6+f_7$. Let $O$ denote the $W$-orbit of $f.$ To see that $\Gamma$ is not $1$-lit, we show that $f_s\not\in O$ for all $s=1,2,\ldots,8.$ Let $\alpha=\alpha_1+\alpha_4+\alpha_5+\alpha_8,$ $\alpha^\vee_1=\alpha_2+\alpha_4+\alpha_5$ and $\alpha^\vee_2=\alpha_1.$ Using (\ref{e4.6}), we find that $\theta(\alpha)=f,$ $\theta(\alpha^\vee_1)=f_1$ and $\theta(\alpha^\vee_2)=f_2.$  Using (\ref{e2.1}) and (\ref{e2.2}), we find that $Q(\alpha)=0$, $Q(\alpha^\vee_1)=1$ and $Q(\alpha^\vee_2)=1.$ By (\ref{e2.7}), neither $\alpha^\vee_1$ nor $\alpha^\vee_2$ is in the $W$-orbit of $\alpha.$ Therefore $f_1\not\in O$ and $f_2\not\in O$ by Corollary~\ref{lem4.4}. By symmetry  $f_i\not\in O$ for $s=3,4,\ldots,8.$

\section{Proof of Theorem~\ref{thm1.4}}

In this section, assume that $\Gamma=(S,R)$ contains at least one edge and fix $x,y\in S$ with $xy\in R$. Define $\widehat{\Gamma}=(\widehat{S},\widehat{R})$ to be the simple graph obtained from $\Gamma$ by inserting a new vertex $z$ on the edge $xy$. In other words, $z$ is an element not in $S$ and the sets $\widehat{S}$ and $\widehat{R}$ are $S\cup\{z\}$ and $R\cup\{xz,yz\}\setminus\{xy\}$, respectively. Let $\widehat{W}$ denote the simply-laced Coxeter group associated with $\widehat{\Gamma}$, namely $\widehat{W}$ is the group generated by all elements of $\widehat{S}$ subject to the following relations:
\begin{align}
s^2&=1,\label{e5.7}\\
(st)^2&=1 \qquad \hbox{if $st\not\in \widehat{R},$} \label{e5.8}\\
(st)^3&=1 \qquad \hbox{if $st\in \widehat{R}$} \label{e5.9}
\end{align}
for all $s,t\in \widehat{S}.$

\begin{lem}\label{thm1.3}
For each $u\in\{x,y\}$ there exists a unique homomorphism $\rho_u:W \rightarrow \widehat{W}$ such that $\rho_u(u)=zuz$ and $\rho_u(s)=s$ for all $s\in S\setminus\{u\}.$
\end{lem}
\begin{proof}
Without loss of generality we assume $u=x.$ We first show the existence of $\rho_x.$ By (\ref{defnW1})--(\ref{defnW3}) it suffices to verify that for all $s,t\in S\setminus\{x\}$,
\begin{align}
s^2&=1, \label{e5.1}\\
(st)^2&=1 \qquad \hbox{if $st\not\in R,$} \label{e5.2}\\
(st)^3&=1 \qquad \hbox{if $st\in R,$} \label{e5.3}\\
(zxz)^2&=1, \label{e5.4}\\
(szxz)^2&=1\qquad  \hbox{if $sx\not\in R,$}  \label{e5.5}\\
(szxz)^3&=1\qquad  \hbox{if $sx\in R$} \label{e5.6}
\end{align}
hold in $\widehat{W}.$ It is clear that (\ref{e5.1})--(\ref{e5.3}) are immediate from (\ref{e5.7})--(\ref{e5.9}), respectively. To obtain (\ref{e5.4}), evaluate the left-hand side of (\ref{e5.4}) using (\ref{e5.7}). By (\ref{e5.8}) and (\ref{e5.9}), for any $s\in S\setminus\{x,y\}$ we have
\begin{align}
(sz)^2&=1, \label{e5.13}\\
(sx)^2&=1 \qquad  \hbox{if $sx\not\in R,$} \label{e5.14}\\
(sx)^3&=1 \qquad  \hbox{if $sx\in R$} \label{e5.15},
\end{align}
and
\begin{align}
(yx)^2&=1, \label{e5.10}\\
(xz)^3&=1,\label{e5.11}\\
(yz)^3&=1 \label{e5.12}
\end{align}
in $\widehat{W}.$ In what follows, the relation (\ref{e5.7}) will henceforth be used tacitly in order to keep the argument concise. Concerning (\ref{e5.5}), let $s\in S\setminus\{x\}$ with $sx\not\in R$ be given. By (\ref{e5.13}), (\ref{e5.14}) the element $s$ commutes with $z$ and $x$ in $\widehat{W},$ respectively. Therefore the left-hand side of (\ref{e5.5}) is equal to $(zxz)^2.$ Now, by (\ref{e5.4}) we have (\ref{e5.5}) in $\widehat W$. To verify (\ref{e5.6}) we divide the argument into the two cases: (A) $s\in S\setminus\{x,y\}$ and $sx\in R;$ (B) $s=y$ in $S.$

\noindent (A) By (\ref{e5.13}) and (\ref{e5.15}), we have $zsz=s$ and $xsxsx=s$ in $\widehat W$, respectively. In the left-hand side of (\ref{e5.6}), replace $zsz$ with $s$ twice and then replace $xsxsx$ with $s.$ This yields $(szxz)^3=(sz)^2$ in $\widehat{W}.$ By (\ref{e5.13}) the relation (\ref{e5.6}) holds.

\noindent (B) In this case we need to show that
\begin{align}\label{e5.17}
(yzxz)^3=1
\end{align}
in $\widehat{W}.$ By (\ref{e5.11}) we have $zxz=xzx$ in $\widehat{W}.$ Use this to rewrite (\ref{e5.17}) as
\begin{align}\label{e5.16}
(yxzx)^3=1.
\end{align}
By (\ref{e5.10}) and (\ref{e5.12}) we have $xyx=y$ and $zyzyz=y$ in $\widehat W$, respectively. In the left-hand side of (\ref{e5.16}), replace $xyx$ with $y$ twice and then replace  $zyzyz$ with $y.$ This yields
$
(yxzx)^3=(yx)^2
$
in $\widehat{W}.$ Now, by (\ref{e5.10}) we have (\ref{e5.16}) in $\widehat W$. Therefore (\ref{e5.6}) holds.

We have shown the existence of $\rho_x$. Such a homomorphism $\rho_x$ is clearly unique since the elements $s\in S$ generate $W$.
\end{proof}

For the rest of this section, let $\rho_x$ and $\rho_y$ be as in Lemma~\ref{thm1.3}.
Let $\widehat{V}$ denote a $\mathbb{F}_2$-vector space that has a basis $\{\alpha_s~|~s\in \widehat{S}\}$ in one-to-one correspondence with $\widehat{S}.$ Let $\widehat{V}^*$ denote the dual space of $\widehat{V}$ and let $\{h_s~|~s\in \widehat{S}\}$ denote the basis of $\widehat{V}^*$ dual to $\{\alpha_s~|~s\in \widehat{S}\}.$ For each $s\in \widehat{S}$ define a linear transformation $\widehat{\kappa}_s:\widehat{V}^*\rightarrow \widehat{V}^*$ by
\begin{align}\label{e6.7}
\widehat{\kappa}_s\hspace{0.5mm}h=h+h(\alpha_s)\sum_{st\in \widehat{R}}h_t \qquad \quad &\hbox{for all $h\in \widehat{V}^*.$}
\end{align}
Let ${\rm GL}(\widehat{V}^*)$ denote the general linear group of $\widehat{V}^*.$ Let $\widehat{\kappa}$ denote the representation from $\widehat{W}$ into ${\rm GL}(\widehat{V}^*)$ such that $\widehat{\kappa}(s)=\widehat{\kappa}_s$ for all $s\in \widehat{S}.$ Define an action of $\widehat{W}$ on $\widehat{V}^*$ by $wh=\widehat{\kappa}(w)h$ for all $w\in \widehat{W}$ and $h\in \widehat{V}^*.$ For each $u\in \{x,y\}$, we define a linear transformation $\delta_u:\widehat{V}^*\to V^*$ by
\begin{align}\label{e6.2}
\delta_u(h_z)=f_u,\qquad \quad \delta_u(h_s)=f_s \qquad \hbox{for all $s\in S.$}
\end{align}
For each $u\in \{x,y\}$ the linear transformation $\delta_u$ is clearly onto and the kernel of $\delta_u$ is
\begin{gather}\label{e6.13}
{\rm Ker}\hspace{0.05cm}\delta_u=\{0,h_u+h_z\}.
\end{gather}
Using (\ref{e6.2}), it is routine to verify that for each $u\in \{x,y\}$ and $s\in S$,
\begin{align}\label{e6.8}
\sum_{st\in \widehat{R}}\delta_u(h_t)=\left\{
\begin{array}{ll}
f_x+f_y+\sum\limits_{ut\in R}f_t \qquad &\hbox{if $s=u,$}\\
\sum\limits_{st\in R}f_t \qquad &\hbox{if $s\not=u.$}
\end{array}
\right.
\end{align}

\begin{lem}\label{lem6.4}
Assume that $O$ is a $\widehat{W}$-orbit of $\widehat{V}^*$ with $O\not=\{0\}.$ Then $\delta_u(O)\not=\{0\}$ for all $u\in\{x,y\}$.
\end{lem}
\begin{proof}
Without loss of generality we show that $\delta_x(O)\not=\{0\}.$  Suppose on the contrary that $\delta_x(O)=\{0\}.$ Since $O\not=\{0\}$ and by (\ref{e6.13}), this forces that $O=\{h_x+h_z\}.$ However $\widehat{\kappa}_z(h_x+h_z)=h_y+h_z\in O$, a contradiction.
\end{proof}

\begin{lem}\label{lem6.3}
For all $u\in\{x,y\}$ and $w\in W$, we have
\begin{align*}
\kappa(w)\circ \delta_u=\delta_u\circ \widehat{\kappa}(\rho_u(w)).
\end{align*}
\end{lem}
\begin{proof}
Let $u\in\{x,y\}$ be given. By Lemma~\ref{thm1.3} and since the elements $s\in S$ generate $W,$ it suffices to show that
\begin{gather}
\kappa_u\circ \delta_u=\delta_u\circ \widehat{\kappa}_z\circ\widehat{\kappa}_u\circ\widehat{\kappa}_z, \label{e6.6}\\
\kappa_s\circ\delta_u=\delta_u\circ\widehat{\kappa}_s \qquad \quad \hbox{for all $s\in S\setminus\{u\}.$} \label{e6.1}
\end{gather}
To verify (\ref{e6.6}), we show that
\begin{align}\label{e6.12}
(\kappa_u\circ \delta_u)(h_s)=(\delta_u\circ \widehat{\kappa}_z\circ\widehat{\kappa}_u\circ\widehat{\kappa}_z)(h_s) \qquad \quad \hbox{for all $s\in \widehat{S}$}.
\end{align}
The argument is divided into the two cases: (A) $s\in\{u,z\};$ (B) $s\in \widehat{S}\setminus\{u,z\}.$

\noindent (A) Using (\ref{e6.7}) we find that $(\widehat{\kappa}_z\circ\widehat{\kappa}_u\circ\widehat{\kappa}_z)(h_s)$ is equal to
$$
h_s+h_x+h_y+\sum\limits_{ut\in \widehat{R}}h_t.
$$
By this and using (\ref{e6.2}) and (\ref{e6.8}), the right-hand side of (\ref{e6.12}) is equal to
\begin{align}\label{e6.5}
f_u+\sum_{ut\in R}f_t.
\end{align}
Using (\ref{e1.3}) and (\ref{e6.2}), the left-hand side of (\ref{e6.12}) is equal to (\ref{e6.5}). Therefore (\ref{e6.12}) holds.

\noindent (B) By (\ref{e1.3}), (\ref{e6.7}) and (\ref{e6.2}), we have $\kappa_u(f_s)=f_s,$ $\kappa_u'(h_s)=\kappa_z'(h_s)=h_s$ and  $\delta_u(h_s)=f_s$, respectively. Using these we find either side of (\ref{e6.12}) is equal to $f_s$. We have shown (\ref{e6.6}).

To verify (\ref{e6.1}), we fix $s\in S\setminus\{u\}$ and show that
\begin{align}\label{e6.11}
(\kappa_s\circ\delta_u) (h_t)=(\delta_u\circ\widehat{\kappa}_s) (h_t) \qquad \quad \hbox{for all $t\in \widehat{S}$}.
\end{align}
The argument is divided into the two cases: (C) $t\in\{u,z\};$ (D) $t\in \widehat{S}\setminus\{u,z\}.$

\noindent (C) By (\ref{e1.3}), (\ref{e6.7}) and (\ref{e6.2}), we have $\kappa_s(f_u)=f_u,$ $\widehat{\kappa}_s(h_t)=h_t$ and $\delta_u(h_t)=f_u$,  respectively. Using these we find either side of (\ref{e6.11}) is equal to $f_u$. Therefore (\ref{e6.11}) holds.

\noindent (D) Using (\ref{e6.7}) and (\ref{e6.8}), the right-hand side of (\ref{e6.11}) is equal to
\begin{align}\label{e6.3}
f_t+h_t(\alpha_s)\sum_{sv\in R}f_v.
\end{align}
Using (\ref{e1.3}) and (\ref{e6.2}), the left-hand side of (\ref{e6.11}) is equal to
\begin{align}\label{e6.4}
f_t+f_t(\alpha_s)\sum_{sv\in R}f_v.
\end{align}
Clearly (\ref{e6.3}) and (\ref{e6.4}) are equal since $f_t(\alpha_s)=h_t(\alpha_s)$. We have shown (\ref{e6.1}). The result follows.
\end{proof}

From now on, assume that $\Gamma=(S,R)$ is a tree with a perfect matching $\mathcal P$. For each $s\in S$, define $A_s$ to be the subset of $S$ consisting of all elements $t\in S\setminus\{s\}$ for which the path in $\Gamma$ joining $s$ and $t$ is an alternating path which starts from and ends on edges in $\mathcal{P}.$ 
Clearly, for each $s\in S$ the number $a_s$ is equal to the size of $A_s$. For each $s\in S$ we define
\begin{align}\label{e7.1}
\alpha^\vee_s=\sum_{t\in A_s}\alpha_t.
\end{align}
Using (\ref{e2.1}), (\ref{e2.2}) and (\ref{e7.1}), the following lemma is straightforward.

\begin{lem}\label{lem7.1}
For each $s\in S$ we have $Q(\alpha_s^\vee)\equiv a_s \pmod 2$.
\end{lem}

\begin{lem}\label{lem7.2}
For each $s\in S$ we have $\theta(\alpha^\vee_s)=f_s$.
\end{lem}
\begin{proof}
Suppose on the contrary that the set $\{A_s~|~s\in S \hbox{ and }\theta(\alpha^\vee_s)\not=f_s\}$ is nonempty. From this set, we choose a minimal element $A_s$ under inclusion. Let $t\in S$ with $st\in \mathcal{P}.$ Observe that $A_s$ is equal to the disjoint union of $\{t\}$ and the sets $A_u$ for all $u\in S\setminus\{s\}$ with $ut\in R$. By this observation and (\ref{e7.1}), we deduce that
\begin{align}\label{e7.2}
\alpha_t=\sum_{ut\in R} \alpha_u^\vee.
\end{align}
By the choice of $A_s$, for each $u\in S\setminus\{s\}$ with $ut\in R$ we have $\theta(\alpha_u^\vee)=f_u$. Apply $\theta$ to either side of (\ref{e7.2}) and then apply (\ref{e4.6}) to the left-hand side. Simplifying the resulting equation, we obtain that $\theta(\alpha_s^\vee)=f_s$, a contradiction.
\end{proof}

We are now ready to prove Theorem~\ref{thm1.4}.

\medskip

\noindent {\bf Proof of Theorem~\ref{thm1.4}:} By Corollary~\ref{prop3.1} the tree $\Gamma$ is nondegenerate. Since all paths are $1$-lit, we may assume that $\Gamma$ is not a path. Let $O$ denote any nonzero $\widehat{W}$-orbit of $\widehat{V}^*$. By Lemma~\ref{lem7.2} we have $\theta(\alpha_x^\vee)=f_x$ and $\theta(\alpha_y^\vee)=f_y$. We first suppose that the edge $xy$ is of odd type. To see that $\widehat\Gamma$ is $1$-lit, it suffices to show that there exists $s\in \widehat{S}$ such that $h_s\in O$. By Lemma~\ref{lem6.4}, there exists $h\in O$ such that $\delta_x(h)\not=0.$ By Lemma~\ref{lem7.1} one of $Q(\alpha_x^\vee)$ and $Q(\alpha_y^\vee)$ is $1$ and the other is $0$. By Corollary~\ref{cor4.1} there exists $w\in W$ such that
$\kappa(w)\delta_x(h)$ is equal to $f_x$ or $f_y$.
By Lemma~\ref{lem6.3} we have
 $\delta_x(\widehat{\kappa}(\rho_x(w))h)$ is equal to $f_x$ or $f_y$.
Using (\ref{e6.13}), we deduce that one of $h_x$, $h_y$, $h_z$, $h_x+h_y+h_z$ is in $O$. By this and since $\widehat{\kappa}_z(h_x+h_y+h_z)=h_z$, one of $h_x$, $h_y$, $h_z$ is in $O$, as desired.

We now suppose that $xy$ is of even type. By Lemma~\ref{lem7.1} we have $Q(\alpha_x^\vee)=Q(\alpha_y^\vee)$. By Corollary~\ref{cor4.1}, the two vectors $f_x$ and $f_y$ are in the same nonzero $W$-orbit of $V^*$ and there exists $u\in S\setminus\{x,y\}$ such that $f_u$ is in the other nonzero $W$-orbit of $V^*$ by Theorem~\ref{thm1.2}. Without loss of generality, we assume that $u$ and $x$ lie in the same component of  the graph $(S,R\setminus\{xy\})$.  By Lemma~\ref{lem6.4}, there exists $h\in O$ such that $\delta_x(h)\not=0.$ By the above comments, there exists $w\in W$ such that $\kappa(w)\delta_x(h)$ is equal to $f_u$ or $f_x$. By Lemma~\ref{lem6.3}, we have $\delta_x(\widehat{\kappa}(\rho_x(w))h)$ is equal to $f_u$ or $f_x$.  Using (\ref{e6.13}), we find that $\widehat{\kappa}(\rho_x(w))h$ is equal to one of $h_u$, $h_x$, $h_z$,  $h_u+h_x+h_z$. In particular $(\widehat{\kappa}(\rho_x(w))h)(\alpha_s)=0$ for all $s\in \widehat{S}\setminus\{u,x,z\}$. If $xy\in \mathcal P$ (resp. $xy\not\in \mathcal P$), we let $\widehat \Gamma_x$ denote the component of the graph $(\widehat S,\widehat R\setminus\{yz\})$ (resp. $(\widehat S,\widehat R\setminus\{xz\})$) containing $x$. Clearly $\widehat \Gamma_x$ is a tree with a perfect matching and contains $u$. Applying Theorem~\ref{thm1.2} to $\widehat \Gamma_x$, there exists a finite sequence of moves for which we only choose the vertices of $\widehat \Gamma_x$ such that $\widehat{\kappa}(\rho_x(w))h$ is transferred to $h'$, where $h'(\alpha_s)=0$ for all $s\in \widehat{S}$ except some vertex in $\widehat \Gamma_x$ and the vertex of $\widehat \Gamma$ that is adjacent to $\widehat \Gamma_x$ and not in $\widehat \Gamma_x$. Therefore $\widehat \Gamma$ is $2$-lit.
\hfill $\Box$

\bigskip

\noindent Hau-wen Huang
\hfil\break Mathematics Division
\hfil\break National Center for Theoretical Sciences
\hfil\break National Tsing-Hua University
\hfil\break Hsinchu 30013, Taiwan, R.O.C.
\hfil\break Email:  {\tt poker80.am94g@nctu.edu.tw}
\medskip

\end{document}